\documentclass[final]{siamltex}
\usepackage{graphics,graphicx}
\usepackage{verbatim}
\def\Re{\hbox{Re}}
\def\C{{\bf C}}
\def\Gc{{G}}
\def\nmax{n_{\max{}}}
\def\pput(#1,#2)#3{\noindent\smash{\raise#2pt\hbox to 0pt
   {\kern #1pt #3\hss}}\ignorespaces}

\title{Solving Laplace problems with
corner singularities via rational functions}

\author{Abinand Gopal\and
Lloyd N.~Trefethen\thanks{\texttt{gopal@maths.ox.ac.uk} and
\texttt{trefethen@maths.ox.ac.uk},
Mathematical Institute, University of Oxford, Oxford, OX2 6GG, UK.}}

\begin{document}

\maketitle

\markboth{\sc Gopal and Trefethen}
{\sc Laplace problems with corner singularities}

\def\sss{\scriptsize}
\def\kk{\kern .7pt {\rm ]}}

\begin{abstract}
A new method is introduced for solving Laplace problems on 2\kern .5pt D
regions with corners by approximation of boundary data by the
real part of a rational function with fixed poles exponentially
clustered near each corner.  Greatly extending a result of
D. J. Newman in 1964 in approximation theory, we first prove that
such approximations can achieve root-exponential convergence for a
wide range of problems, all the way up to the corner singularities.
We then develop a numerical method to compute approximations via
linear least-squares fitting on the boundary.  Typical problems are
solved in ${<}\kern 1pt 1$s on a desktop to 8-digit accuracy, with
the accuracy guaranteed in the interior by the maximum principle.  The computed
solution is represented globally by a single formula, which can
be evaluated in a few microseconds at each point.

\end{abstract}



\pagestyle{myheadings}
\thispagestyle{plain}

\section{Introduction}
In 1964 Donald Newman published a result that has attracted a good
deal of attention among approximation theorists~\cite{newman}.
Newman considered supremum norm approximation of $f(x) = |x|$
on $[-1,1]$ and showed that, whereas degree $n$ polynomial
approximants offer at best $O(n^{-1})$ convergence as $n\to\infty$,
rational functions can achieve {\em root-exponential convergence,}
that is, errors $O(\exp(-C\sqrt n\kern 1.2pt))$ with $C>0$.
(The degree of a rational function is the maximum of the degrees
of its numerator and denominator.)  Newman's construction
crucially involved placing poles in the complex plane clustered
exponentially near the singularity at $x=0$, and he also
showed that root-exponential convergence was the best possible.
In the subsequent half-century, his result has been sharpened
by approximation theorists, notably Herbert Stahl, and extended
to other functions such as $|x|^\alpha$~\cite{gon67,stahl,stenger}.
However, there appear to have been no attempts heretofore to
exploit the root-exponential effect in scientific computing.

Here we introduce a method of this kind for solving Laplace
problems on a 2\kern .5pt D domain $\Omega$ bounded by 
a piecewise smooth Jordan curve $\Gamma$ with corners,
such as a polygon.
In the simplest situation we have a Dirichlet problem with
specified real boundary data $h$,
\begin{equation}
\Delta u (z) = 0, ~ z\in \Omega,
\qquad u(z) = h(z), ~ z\in \Gamma,
\label{theproblem}
\end{equation}
where we use complex notation $z = x+iy$ for simplicity.
The method approximates $u$ by the real part of a rational function,
\begin{equation}
u(z) \approx \Re \kern 1.2pt r(z),
\label{realpart}
\end{equation}
with $r$ taking the form
\begin{equation}
r(z) \;=\;\; \sum_{j=1}^{N_1}{a_j\over z-z_j}
\pput(-36,22){\scriptsize\sl ``\kern 1pt NEWMAN\kern 1pt''}
\;\;+\; \sum_{j=0}^{N_2}b_j (z-z_*^{})^j.
\pput(-48,22){\scriptsize\sl ``\kern 1pt RUNGE\kern 1pt''}
\label{form}
\end{equation}
Computing an optimal rational function $r$ is a difficult nonlinear
problem in general~\cite{hochman,aaa,atap}, but our method is
linear because it fixes the poles $\{z_j\}$ of the first sum
in (\ref{form}) a priori in a configuration with exponential
clustering near each corner.  We speak of this as the ``Newman''
part of (\ref{form}), concerned with resolving singularities.
The second sum is simply a polynomial, with the expansion point
$z_*$ taken as a point roughly in the middle of $\kern 1pt\Omega$,
and we call this the ``Runge'' term, targeting the smoother part
of the problem.  (The observation that functions analytic in a
neighborhood of a simply connected compact set in the complex plane
can be approximated by polynomials with exponential convergence
dates back to Runge in 1885~\cite{gaier,runge,series}.)  With $\{z_j\}$
and $z_*^{}$ fixed, our method finds coefficients $a_j$ and $b_j$
by solving a linear least-squares problem on a discrete subset
of points on the boundary~$\Gamma$, which are also exponentially
clustered in a manner reflecting that of $\{z_j\}$.

One might imagine that achieving root-exponential convergence would
require a delicate choice of poles $\{z_j\}$ in (\ref{form}).
However, virtually any exponential clustering is in fact
sufficient, provided it scales with $n^{-1/2}$ as $n\to\infty$,
and Section~2 is devoted to presenting theorems to establish
this claim.  Quite apart from their application to Laplace
problems, we believe these results represent a significant
addition to the approximation theory literature, as well as
shedding light on the clustered poles observed experimentally
in~\cite{conf} and~\cite{hochman}.  Section~3 describes our
algorithm, which depends on placing sample points on the boundary
with exponential clustering to match that of the poles outside.
Section~4 presents examples illustrating its remarkable speed
and accuracy, and Section~5 comments on variants such as Neumann
boundary conditions and multiply connected domains.  We are in
the process of developing codes to make
the method available to others, and an analogous method for the
Helmholtz equation will be presented in a separate publication.

A brief announcement of the new method, without theorems or
proofs, was published in~\cite{pnas}.  Our method is a variant
of the Method of Fundamental Solutions (MFS), and a discussion
of its relationship with MFS and other methods is given
in the final section.

\section{\label{thms}Root-exponential convergence theorems}
This section will show that root-exponential convergence is
achievable for all kinds of domains provided that poles are
clustered exponentially near corner singularities with a spacing
that scales as $O(n^{-1/2})$.  The result is sharp in the sense
that faster than root-exponential convergence is in general not
possible for rational approximations of Laplace problems with
corner singularities.  To prove this, it is enough to consider
the problem (\ref{theproblem}) with $\Omega$ as the upper half
of the unit disk and $h(z) = |\Re z|$.  If $\Re \kern 1.5pt
r(z) \approx u(z)$ for $z\in\Omega$ for some rational function
$r$, then in particular $\Re\kern 1.5pt r(x)\approx |x|$ for
$x\in [-1,1]$.  But for $x\in [-1,1]$, $\Re\kern 1.5pt r(x) =
(r(z)+\overline{r(z)}\kern .7pt )/2$ is itself a rational function
of $x$.  Thus by Newman's converse result, the approximation can
improve no faster than root-exponentially with $n$.

Our proofs combine the Cauchy integral formula, to decompose
global problems into local pieces, and the Hermite integral formula
for rational interpolation.  The rational version of the Hermite
formula, which goes back at least to J. L. Walsh in the 1930\kern
.5pt s~\cite[sec.~8.1]{walsh}, is well known by approximation
theorists, and it is possible that results like ours can be found
in the literature.  Examples of related works are~\cite{bagbygaut}
and~\cite{stenger}.  However, the emphasis in the literature is
on relatively delicate analysis in relatively specific settings.
(This comment applies to Theorem~1 of our own paper~\cite{conf},
where we used an exponential change of variables coupled with the
equispaced trapezoidal rule to prove root-exponential convergence
for approximations of $x^\alpha$ on a half-disk.)  Here we will
make more general arguments.  The root-exponential rate is a
consequence of the balance of local accuracy at a corner, which
requires poles to be placed closely nearby, and global accuracy
further from the corner, which requires poles to be spaced close
together.  This is the same balance that leads to root-exponential
convergence of exponential quadrature formulas~\cite{mori,trap}.

We begin with a statement of the Hermite formula.  A rational
function $r$ is said to be of {\em type $(m,n)$} if it can be
written in the form $r=p/q$ where $p$ and $q$ are polynomials of
degrees $\le m$ and $\le n$, respectively.

\smallskip

\begin{lemma}
\underline{Hermite integral formula for rational interpolation}.
Let\/ $\Omega$ be a simply connected domain in $\C$ bounded
by a closed curve $\Gamma$, and let $f$ be analytic in 
$\Omega$ and extend continuously to the boundary.
Let interpolation points $\alpha_0,\dots, \alpha_{n-1}\in \Omega$
and poles $\beta_0,\dots,\beta_{n-1}$ anywhere in
the complex plane be given.
Let $r$ be the unique type $(n-1,n)$ rational function
with simple poles at $\{\beta_j\}$ that
interpolates $f$ at $\{\alpha_j\}$.  Then for any $z\in\Omega$,
\begin{equation}
f(z) - r(z) = {1\over 2\pi i} \int_\Gamma {\phi(z)\over \phi(t)}
{f(t)\over t-z} \kern 1pt dt,
\label{herm}
\end{equation}
where
\begin{equation}
\phi(z) = \prod_{j=0}^{n-1}(z-\alpha_j)\biggl/\kern 1.5pt
\prod_{j=0}^{n-1}(z-\beta_j).
\label{phidef}
\end{equation}
\label{hermlemma}
\end{lemma}

\smallskip

\begin{proof}
The theorem is presented as Theorem 2 of Chapter 8 of
Walsh's monograph~\cite{walsh}, with a proof based
on residue calculus.  Walsh quoted the same formula in
an earlier paper of his as equation (12.1) of~\cite{walsh32}.
We do not know whether he was the first to develop this result.
For a more recent discussion, see~\cite{levsaff}.
\end{proof}
\smallskip

The application of (\ref{herm}) for accuracy estimates goes back
to the analysis of polynomial interpolants by M\'eray and Runge
at the end of the 19th century~\cite[chap.~11]{atap}, and in
the case of rational functions, to Walsh~\cite{walsh32, walsh}.  With good
choices of $\{\alpha_j\}$ and $\{\beta_j\}$, $\phi$ may be much
larger for $t\in\Gamma$ than for $z\in\Omega$.  In such a case,
the ratio $\phi(z)/\phi(t)$ in (\ref{herm}) will be very small,
and this enables one to use the integral to bound $f(z) - r(z)$.
For the estimate of our first theorem, we will apply just this
reasoning.

Note that unlike the proofs of our theorems, our algorithms are
not based on interpolation.  The purpose of the theorems is to
establish the existence of root-exponentially good approximants,
not to construct them.  That will be done by the much more flexible
method of least-squares fitting on the boundary.

\begin{figure}
\centering
\vskip .3in
\includegraphics[scale=.9]{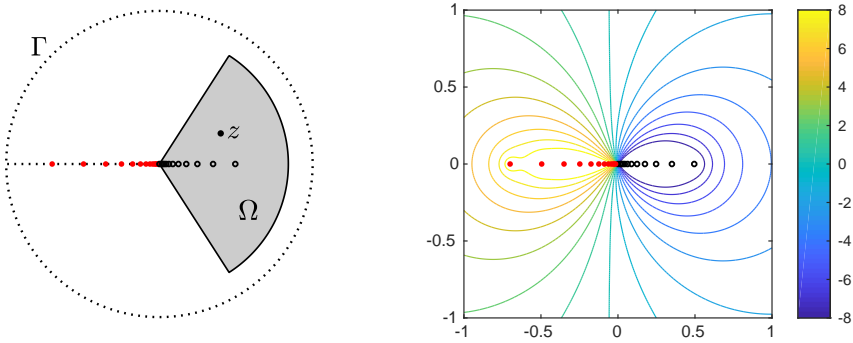}
\caption{On the left, the approximation problem of
Theorem~$\ref{thm1}$.  A bounded analytic function~$f$ is given in
the unit disk slit along the negative
real axis.  Root-exponentially accurate
rational approximants to $f$ in a wedge
$\Omega = \rho A_\theta$ are constructed
by placing poles in $[-1,0)$ exponentially
clustered near $x=0$
and interpolation points similarly clustered in $[\kern .3pt 0,1)$.
The image on the right shows level curves of $\log_{10}|\phi(z)|$,
where $\phi$ is defined by $(\ref{phidef})$,
$(\ref{poles})$, and $(\ref{interppts})$ with $n=12$, $\sigma = 0.1$.}
\label{fig1}
\end{figure}

We begin with a rational approximation problem on the slice-of-pie
region sketched on the left in Figure~\ref{fig1}, 
starting from the definition
\begin{equation}
A_\theta = \{z\in\C:\, |z|<1,\,-\theta< \hbox{\rm arg\kern 1pt} z <
\theta\},
\end{equation}
and assume a function $f$ is given that is bounded and analytic in
the slit disk $A_\pi$.  For example, $f$ might have a $z^\sigma$
or $z^\sigma \log z$ singularity at $z=0$ ($\sigma>0$).
As shown in~\cite{lehman} and \cite{wasow}, these are the
standard singularities that arise in Laplace problems bounded
by analytic curves meeting at corners.  The following theorem,
which is close to Theorem~2.1a of~\cite{stenger}, approximates
on a wedge with half-angle $\theta <\kern 1pt \pi/2$.  However,
the result actually holds for any $\theta<\pi$, as we discuss
at the end of the section.  The notation $\|\cdot\|_\Omega^{}$
denotes the supremum norm over $\Omega$.

\smallskip

\begin{theorem}
\label{thm1}
Let $f$ be a bounded analytic function in the slit disk $A_\pi$
that satisfies $f(z) = O(|z|^\delta)$ as $z\to 0$ for some $\delta>0$,
and let $\theta\in (0,\pi/2)$ be fixed.
Then for some $\rho\in (0,1)$ depending on\/ $\theta$ but not $f$, there
exist type $(n-1,n)$ rational
functions $\{r_n\},$ $1\le n < \infty$, such that
\begin{equation}
\|f-r_n\|_\Omega^{} = O(e^{-C\sqrt n \kern 1pt})
\label{mainbound}
\end{equation}
as $n\to\infty$ for some $C>0$, where $\Omega = \rho A_\theta$.
Moreover, each\/ $r_n$ can be taken to have simple poles only at
\begin{equation}
\beta_j = -e^{-\sigma j/\sqrt n\kern 1.5pt}, \quad 0\le j \le n-1,
\label{poles}
\end{equation}
where $\sigma >0$ is arbitrary.
\end{theorem}

\smallskip

\begin{proof}
Fix $\sigma >0$ and define
poles $\beta_j$ by (\ref{poles}) and interpolation points by
\begin{equation}
	\alpha_0^{} = 0, \qquad \alpha_j = -\beta_j, ~ 1\le j \le n-1 .
\label{interppts}
\end{equation}
This implies
$|z-\alpha_0^{}|/|z-\beta_0^{}| \le |z|$ and
$|z-\alpha_j^{}|/|z-\beta_j^{}| \le 1$ for $j\ge 1$ for
any $z\in\Omega$, and therefore by (\ref{phidef})
\begin{equation}
|\phi(z)| \le |z|, \quad z\in\Omega.
\label{pairs1}
\end{equation}
Here and until the final paragraph of the proof, the constant
$\rho\in(0,1)$ has not yet been chosen, and our statements about
$\Omega = \rho\kern .5pt A_\theta$ apply regardless of its value.

Let $r$ be the type $(n-1,n)$
rational interpolant to $f$ as in Lemma~\ref{hermlemma}.
We need to show that there are constants $A,C>0$, independent
of $z$, such that
\begin{equation}
|f(z) - r(z)| \le A e^{-C\sqrt n } , \quad n\ge 1
\label{rootexpbound}
\end{equation}
for $z\in \Omega$.  To do this, we apply the integral (\ref{herm}),
with $\Gamma$ as the boundary of the slit disk.  (The two sides of
$[-1,0\kern .3pt ]$ are distinct components of $\Gamma$.  We may
assume without loss of generality that $f$ extends continuously
to the boundary, for if not, we may integrate over contours in
$A_\pi$ that come arbitrarily close to $\Gamma$.)

\begin{figure}
\begin{center}
\vskip .15in
\includegraphics{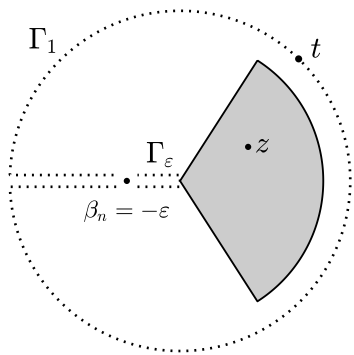}
\end{center}
\caption{\label{sketch}Decomposition of\/ $\Gamma$
into two pieces to prove by\/ $(\ref{hermsplit})$ that
$f(z)-r(z)$ is root-exponentially small for
$z\in\Omega$.  On $\Gamma_\varepsilon$, the integral
is small because $f$ is small and the contour is short.
On $\Gamma_1$, it is
small because $\phi(z)/\phi(t)$ is small.}
\end{figure}

Now for each $n$, by (\ref{poles}),
the pole closest to zero is
$\beta_n= -\varepsilon$ with $\varepsilon = \exp(-\sigma (n-1)/\sqrt n\kern 1.5pt)$.
We split $\Gamma$ into two pieces,
as sketched in Figure~\ref{sketch}:
\begin{displaymath}
\Gamma_\varepsilon = \{t\in \Gamma: \; |t|<\varepsilon\},
\quad
\Gamma_1 = \Gamma\backslash \Gamma_\varepsilon.
\end{displaymath}
Equation (\ref{herm}) becomes
\begin{equation}
f(z) - r(z) = I_\varepsilon + I_1 =
{1\over 2\pi i} \int_{\Gamma_\varepsilon} {\phi(z)\over \phi(t)}
{f(t)\over t-z} \kern 1pt dt
+ {1\over 2\pi i} \int_{\Gamma_1} {\phi(z)\over \phi(t)}
{f(t)\over t-z} \kern 1pt dt,
\label{hermsplit}
\end{equation}
and we need to show that $|I_\varepsilon|$ and $|I_1|$ are both bounded
by $Ae^{-C\sqrt n}$.

To bound $I_\varepsilon$, we note that for
$t\in \Gamma_\varepsilon$, we have $|z|\le |\kern .5pt t-z|$, and thus 
by (\ref{pairs1}),
\begin{displaymath}
\left|{\phi(z)\over t-z}\right| \le 1.
\end{displaymath}
By the H\"older inequality, this
means that to bound $I_\varepsilon$, it is enough to bound
the integral of $f(t)/\phi(t)$ over $\Gamma_\varepsilon$.
Now in analogy to (\ref{pairs1}) we have
\begin{equation}
|\phi(t)| \ge |\kern .3pt t|
\label{pairs}
\end{equation}
for $t\in\Gamma_\varepsilon$, since $|t-\alpha_0|/|t-\beta_0|\ge |t|$
and
$|t-\alpha_j|/|t-\beta_j|\ge 1$ for $j\ge 1$ in (\ref{phidef}).
So it is enough to bound the integral of $f(t)/t$, and by the assumption
on the behavior of $f(z)$ as $z\to 0$,
this amounts to an integral of the type
$$
\int_0^\varepsilon {t^\delta\over t} \kern 1pt dt = 
{\varepsilon^\delta\over \delta}.
$$
Since $\varepsilon = \exp(-\sigma (n-1)/\sqrt n\kern 1.5pt)$, this is
of order $\exp(-\sigma \kern .3pt \delta \sqrt n\kern 1.5pt)$,
independently of $z$ as required.
\vskip 2pt

To bound $I_1$ in (\ref{hermsplit}), we note that although
$f(t)/(t-z)$ grows large pointwise on $\Gamma_1$ for $t\approx
-\varepsilon$ and $z\approx 0$, its integral over $\Gamma_1$ is
bounded, independently of $z\in\Omega$.  Thus by the H\"older
inequality again, it is enough to show that $\phi(z)/\phi(t)$
is root-exponentially small uniformly for $z\in\Omega$ and
$t\in\Gamma_1$.  To do this, we first observe that $\phi(z)$ is
root-exponentially small for $z\in\Omega$.  For $|z|<2\kern .3pt
\varepsilon$, this conclusion follows from (\ref{pairs1}).  For
$|z|\ge 2\kern .3pt \varepsilon$, we note from (\ref{phidef}) that
$\phi(z)$ is a product of factors $(z-\alpha_j)/(z-\beta_j)$ of
size at most $1$.  To show that the product is root-exponentially
small, it suffices to show that on the order of $\sqrt n$
of these factors are bounded in absolute value by a fixed
constant $D<1$ independent of $z$.  For each $z$, a suitable set
of factors with this property are those with $|z|/2 < \alpha_j
< |z|$.  From (\ref{poles}) it follows that the number of these
factors grows in proportion to $\sqrt n$ as $n\to\infty$, and from
elementary geometry it follows that each factor is bounded by a
constant $D<1$ (dependent on $\theta$ but not $z$) as required.

We then have to balance the size of $\phi(z)$ against that
of $\phi(t)$.  For $t\in\Gamma_1$ in the left half-plane,
this is immediate since (\ref{pairs}) ensures that $\phi(t)$
is not small (more precisely, (\ref{pairs}) weakens slightly to
$|\phi(t)|\ge |\kern .3pt t|/\sqrt 2$ for all $t\in \Gamma_1$ in
the left half-plane).  For $t\in\Gamma_1$ in the right half-plane,
we have to be more careful since $\phi(t)$ is root-exponentially
small like $\phi(z)$, and this is where the constant $\rho\in(0,1)$
comes into play. Given $\rho$, let us discard all poles with
$|\beta_j|>\rho$ and interpolation points with $|\alpha_j|>\rho$
(both for $\kern .7pt j\ge 1$) from the rational interpolation problem.  The
remaining quotients $(t-\alpha_j)/(t-\beta_j)$ defining $\phi(t)$
are each of absolute value no smaller than $1-O(\kern .7pt\rho)$,
and their product is no smaller than $(1-O(\kern .7pt\rho))^{\sqrt
n}$, the exponent here being $\sqrt n$ rather than $n$ because the
ratios $(t-\alpha_j)/(t-\beta_j)$ approach $1$ geometrically as a
function of $j/\sqrt n$, so that, loosely speaking, only the first
$O(\sqrt n\kern 1pt)$ of them affect the size of their product.
This implies that the root-exponential constant of $\phi(t)$ in
(\ref{rootexpbound}) satisfies $C = O(\kern .7pt\rho)$.  So for
sufficiently small $\rho$, the rate of root-exponential decrease
of $\phi(z)$ as $n\to\infty$ exceeds that of $\phi(t)$, making
$\phi(z)/\phi(t)$ root-exponentially small.  \end{proof}

\smallskip

\begin{figure}
\centering
\vskip .3in
\includegraphics[scale=1.1]{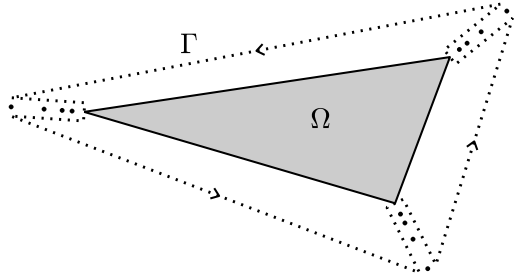}
\caption{The approximation problem of Theorem~$\ref{thm2}$.  An analytic
function\/ $f$ on a convex \hbox{$m$-gon} \kern 1pt $\Omega$ is decomposed as the sum
of\/ $2m$ Cauchy integrals:\ one along the two sides of a slit exterior to
each corner, and one along each line segment connecting the ends of those
slit contours.  
The\/ $m$ corner functions are approximated by
rational functions following Theorem~$\ref{thm1}$, and the\/ $m$ connecting
functions are approximated by a polynomial.}
\label{fig2}
\end{figure}

Theorem~\ref{thm1} is the local assertion, establishing
root-exponential resolution of corner singularities.  The next
step is to show that global approximations can be constructed
by adding together these local pieces, with a polynomial term
to handle smooth components away from the corners.  We make
the argument for the case in which $\Omega$ is a convex polygon
(with internal angles $<\pi$), as sketched in Figure~\ref{fig2},
although in fact we believe convexity is not necessary.

\smallskip
\begin{theorem}
\label{thm2}
Let\/ $\Omega$ be a convex polygon with corners $w_1,\dots, w_m$,
and let
$f$ be an analytic function in $\Omega$ that is analytic on the
interior of each side segment and can be analytically continued to 
a disk near each\/ $w_k$ with a slit along the exterior bisector there.
Assume $f$ satisfies $f(z) - f(w_k) = 
O(|z-w_k|^\delta)$ as $z\to w_k$ for each~$k$ for some $\delta>0$.
There exist degree $n$ rational 
functions $\{r_n\},$ $1\le n < \infty$, such that
\begin{equation}
\|f-r_n\|_\Omega^{} = O(e^{-C\sqrt n \kern 1pt})
\label{mainbound2}
\end{equation}
as $n\to\infty$ for some $C>0$.  Moreover, each $r_n$ can be taken
to have finite poles only at points exponentially clustered along
the exterior bisectors at the corners, with arbitrary clustering
parameter $\sigma$ as in $(\ref{poles})$, so long as the number
of poles near each $w_k$ grows at least in proportion to $n$
as $n\to\infty$.  \end{theorem}

\smallskip
\begin{proof}
We represent $f$ as a sum of $2m$ terms,
\begin{equation}
f = \sum_{k=1}^m f_k + \sum_{k=1}^m g_k,
\label{pieces}
\end{equation}
each defined by a Cauchy integral
\begin{equation}
{1\over 2\pi i} \int_\Gc {f(t)\over t-z} \kern 1pt dt
\label{cauchyint}
\end{equation}
over a different contour $G$ within the closure of the
region of analyticity
of $f$.  For $f_k$, $\Gc$ consists of the two sides of an exterior
bisector at $w_k$, and for $g_k$, $\Gc$ is a curve that
connects the end of the slit contour at vertex $k$ to the beginning
of the slit contour at vertex $k+1$.  See Figure~\ref{fig2}.
The decomposition (\ref{pieces}) holds since the sum of the
contours winds all the way around $\Omega$ while remaining in
the region of analyticity of $f$.  For the mathematics of Cauchy
integrals over open arcs, see Chapter 14 of~\cite{acca}.

Consider first a function $g_k$ defined by the Cauchy integral
(\ref{cauchyint}) over an arc $G$ connecting two adjacent slits.
This function is analytic in any open region disjoint from $G$,
and in particular, in an open region containing the closure
of $\Omega$.  By Runge's theorem~\cite{gaier,runge}, it follows
that $g_k$ can be approximated by polynomials on $\Omega$ with
exponential convergence.  This is the smooth part of the problem.

For the ``Newman'' part, consider a function $f_k$ defined by the
integral over a contour $G$ consisting of two sides of a slit
at $w_k$ associated with a sector of half-angle $\theta<\pi $.
Again, $f_k$ is analytic in any open region disjoint from $G$.
In particular, it is analytic in a slit-disk region like
$A_\pi$ of Figures~\ref{fig1} and~\ref{sketch},
but tilted and translated to lie around the slit $G$ and of
arbitrary scale $R>0$.  If~$\rho$ is the parameter in the proof of
Theorem~\ref{thm1} for the given angle $\theta$, then taking $R$
to be greater than the diameter of $\Omega$ divided by $\rho$
is enough for the argument of Theorem~\ref{thm1} to ensure
the existence of rational approximations to $f_k$ converging
root-exponentially on $\Omega$.  We note that this argument
requires $f_k(z) = O(|z-w_k|^\delta)$ as $z\to w_k$.  This we
can ensure by subtracting a constant from $f_k$, to be absorbed
in the polynomial term of the approximation (\ref{pieces}).
The necessary continuity condition on $f_k$ then follows from
the continuity assumption on $f$ since $f$ and $f_k$ differ by
terms that are analytic at $w_k$.

Adding the $2\kern .3pt m$ polynomial and rational pieces together,
we get rational approximations to $f$ with root-exponential 
convergence in $\Omega$.
\end{proof}
\smallskip

With Theorem~\ref{thm2} we have established the existence of
root-exponentially convergent rational approximations to an
analytic function $f$ on a convex polygonal region $\Omega$,
provided the singularities of $f$ on the boundary are just branch
points at the corners.  We now discuss three extensions.

First of all, nothing in our arguments has utilized the
straightness of the sides of $\Omega$.  The same results hold if
$\Omega$ is bounded by analytic arcs meeting at corners, which is
the setting of the paper by Wasow mentioned earlier~\cite{wasow}
as well as of related papers by Lewy~\cite{lewy59} and
Lehman~\cite{lehman}, who generalized from the Laplace equation
to other PDE\kern .3pt s.  The case of analytic arcs meeting with
angle zero, a cusp, is also not a problem, and indeed in our
context is indistinguishable from the case of a finite angle,
in view of the assumption of analyticity in a slit disk around
each corner.

Secondly, our application of these results to solve the
Laplace equation will concern harmonic functions rather than
analytic ones.  There is little difference in the two settings,
however, since any harmonic function $u$ on $\Omega$ is the real
part of an analytic function $f = u + iv$ on $\Omega$ that is
uniquely determined up to an additive constant (since $\Omega$
is simply connected).  The only technical issue to be considered
is that of regularity.  If $u$ is harmonic in $\Omega$ and can
be harmonically extended across the boundaries and around slits
at each corner, with $u(z)- u(w_k)$ = $O((z-w_k)^\delta)$ for
some $\delta>0$ for each~$k$, does this imply the same properties
for $v$ and thus for the function $f$ to which we would like to
apply Theorem~\ref{thm2}\kern 1pt?  The answer is yes, as can be
shown by a conformal transplantation from $\Omega$ to a disk or
a half-plane following by application of standard theory of
the Hilbert transform, which maps a function to its harmonic
conjugate.  See Theorem~14.2c of~\cite{acca}.

The third issue to be considered is non-convex
domains.  In fact, we believe Theorem~\ref{thm1} is valid without
the assumption $\theta<\pi/2$ and Theorem~\ref{thm2} is valid without
the assumption of convexity.  We will not attempt
a proof, but will outline the very interesting mathematics that
could be applied to this end.

\begin{figure}
\centering
\includegraphics[scale=.9]{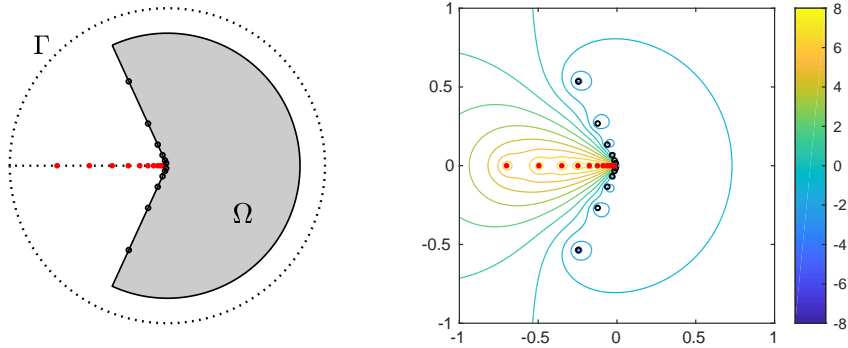}
\caption{For a non-convex sector, one can cluster
interpolation points exponentially on both sides.
More complicated domains than this can be handled
by approximately minimizing the energy\/~$(\ref{energy})$.}
\label{fig3}
\end{figure}

The reason why Theorem~\ref{thm1} is restricted to convex corners
is that its proof makes use of interpolation points $\{\alpha_j\}$
placed along the interior bisector $[\kern .5pt 0,1)$, a convenient
configuration for a simple argument.  To handle more general cases,
better choices of interpolation points are needed.  One interesting
choice is to place points with exponential clustering not on the
bisector, but on the two sides of the sector.  By this means one
can extend Theorem~\ref{thm1} to all $\theta<\pi$, as illustrated
in Figure~\ref{fig3}.  That is,
one can guarantee the existence of rational approximations
that converge root-exponentially near the singular corner of a
non-convex sector of a disk.

For a generalization of Theorem~\ref{thm2} to arbitrary nonconvex
polygons or their curved analogues, however, more is
needed.  An arbitrary region $\Omega$ of this type might have a complicated
shape like an S or a spiral, and to handle such cases, we need to
construct rational approximations of each corner function $f_k$ that
are accurate globally in $\Omega$, not just in a sector.  
To achieve this, the idea we need is to let the interpolation
points $\{\alpha_j\}$ distribute themselves in a configuration that
is approximately optimal in the sense of minimizing 
a certain electrostatic energy associated
with the poles and interpolation points.  Such ideas came to the
fore in approximation theory with the work of Gonchar and
his collaborators beginning shortly after Newman's paper
appeared~\cite{gon67,gonchar,levsaff,safftotik,stahl}.
Given any choice of $\{\alpha_j\}$ and $\{\beta_j\}$,
define the associated {\em potential function} by
\begin{equation}
\Phi(z) \,=\, \sum_{j=0}^{n-1} \log |z-\alpha_j| -
\sum_{j=0}^{n-1} \log|z-\beta_j| ,
\label{potential}
\end{equation}
which is related to the function $\phi$ of (\ref{herm}) by
\begin{equation}
\Phi(z) = \log |\phi(z)|.
\label{hermpphi}
\end{equation}
The point of Theorem~\ref{thm1} is that with good choices
of $\{\alpha_j\}$ and $\{\beta_j\}$, we can ensure that $\Phi$
is smaller on $\Omega$ than on $\Gamma$, so that (\ref{herm})
implies that rational interpolation gives accurate approximants.

In rational approximation theory~\cite{levsaff}, it is usual
to let both $\{\alpha_j\}$ and $\{\beta_j\}$ vary
to find near-optimal approximants by minimizing the energy
\begin{displaymath}
I(\alpha,\beta) \,=\, \sum_{j=0}^{n-1}\kern 3pt \sum_{k=0}^{n-1}
\log|\alpha_j-\beta_k|
-\sum_{j=0}^{n-2}\kern 0pt \sum_{k=j+1}^{n-1}
\Big(\kern -2pt \log|\alpha_j-\alpha_k|+\log|\kern .5pt \beta_j-\beta_k|\Big).
\end{displaymath}
For our purposes, however, the poles $\{\beta_j\}$ are fixed
a priori in a configuration with exponential clustering
near the corners.   The terms involving $\log|\kern .5pt
\beta_j-\beta_k|$ accordingly contribute just a constant,
and the issue reduces to the
adjustment of the interpolation points $\{\alpha_j\}$ to
approximately minimize the energy
\begin{equation}
I(\alpha) \,=\, \sum_{j=0}^{n-1}\kern 3pt \sum_{k=0}^{n-1}
\log|\alpha_j-\beta_k|
-\sum_{j=0}^{n-2}\kern 0pt \sum_{k=j+1}^{n-1}
\log|\alpha_j-\alpha_k|.
\label{energy}
\end{equation} By means of such a minimization, rational
approximations can be constructed that are good throughout a
nonconvex region $\Omega$ with corners.  The convergence rate will
be root-exponential, though with a suboptimal constant since the
poles have not been chosen optimally.

At the end of seven pages of perhaps dense mathematics we remind
the reader that none of these details of Hermite integrals,
Cauchy integrals, and potential theory are part of our algorithm,
which works simply by solving a linear least-squares problem.
Their only purpose is to justify the algorithm theoretically.

\section{The algorithm}
Here is our algorithm for solving (\ref{theproblem}).
At present our codes are exploratory, but we aim to develop more
robustly engineered software in the near future.  \medskip

{\parindent=20pt\parskip=1pt
\def\\{\indent\kern .2in}\def\kk{\kern .5pt ;}
\obeylines\em
ALGORITHM
\kern 1pt
1. Define boundary\/ $\Gamma$, corners $w_1,\dots ,w_m$, boundary function $h$, and \rlap{tolerance $\varepsilon$.}
\kern 1pt
2. For increasing values of\/ $n$ with $\sqrt n$ approximately evenly spaced\kern .5pt:
\\ 2a. Fix $N_1 = O(mn)$ poles $1/(z-z_k)$ clustered outside the corners\kk
\\ 2b. Fix $N_2+1 = O(n)$ monomials $1,\dots, (z-z_*)^{N_2}$ and set $N=2N_1+2N_2+1$;
\\ 2c. Choose $M\approx 3N$ sample points on boundary, also clustered near corners\kk
\\ 2d. Evaluate at sample points to obtain $M\times N$ matrix $A$ and RHS vector $b$\kk
\\ 2e. Solve the least-squares problem $Ax\approx b$ for the coefficient vector $x$\kk
\\ 2f. Exit loop if\/ $\|Ax-b\|_\infty < \varepsilon$ or if\/ $N$ is too large or the error is growing.
\kern 1pt
3. Confirm accuracy by checking the error on a finer boundary mesh.
\kern 1pt
4. Construct a function to evaluate $r(z)$ based on the computed coefficients $x$.
\par}

\medskip

If $\nmax$ is the largest value of\/ $n$ used in the computation, then the
largest matrix will be of row and column dimensions $O(m\kern .5pt \nmax)$,
corresponding to an operation count in step
(2e) of $O((m\kern .5pt \nmax)^3)$.  If we
assume $\nmax = O(|\log \varepsilon|^2)$, this gives an overall operation
count for the algorithm of
\begin{equation}
\hbox{Operation count} = O(m^3 |\log(\varepsilon)|^6).
\label{opcount}
\end{equation}
These exponents look daunting, and we are investigating methods
of speeding up the linear algebra.  (Both domain decomposition
by corners and multiscale decomposition by space scales offer
possibilities.)  However, the constants in the ``$O$'' are
small, and small- or medium-scale problems often fall short
of the asymptotic regime of cubic linear algebra in any case.
As we shall show in the next section, our MATLAB implementation
typically solves problems with $m \le 8$ and $\varepsilon\ge
10^{-8}$ in a fraction of a second on a desktop machine, with the accuracy
guaranteed all the way up to the corners.

We now give details of each of the steps.

\medskip
{\em 1. Define boundary $\Gamma$, corners $w_1,\dots ,w_m$,
boundary function $h$, and tolerance $\varepsilon$.}\hfill\break
If \kern 1.5pt $\Gamma$ is a polygon, these setup operations
are straightforward, and for more complicated domains,
any convenient representation of the boundary can be used.
Concerning~$\varepsilon$, our usual habit is to work by default
with the value $10^{-6}$, which leads to solutions very often
in less than a tenth of a second.  The estimate (\ref{opcount}) would
suggest that tightening $\varepsilon$ to $10^{-8}$ should make a
computation ${\approx}\kern 1pt 5.6$ times slower, but in practice,
for problems of modest scale, the slowdown is usually more like
a factor of $2$.

{\em 2a. Fix\/ $N_1 = O(mn)$ poles $1/(z-z_k)$ clustered outside
the corners.}
In all the experiments of this paper, following a rule of
thumb derived from numerical experience, $n$ poles are placed
near each salient corner and $3n$ poles 
at each reentrant corner (i.e., with internal
angle $>\pi$).  
In the version of the paper originally submitted for publication,
the poles were clustered according to (\ref{poles}) with $\sigma = 2.5$.
Subsequent numerical experience showed that the numbers of degrees
of freedom could be reduced about $20\%$ by modifying the prescription 
to
\begin{equation}
\beta_j = -e^{-\sigma (\sqrt n-\sqrt j\kern 1pt )}, \quad 1\le j \le n
\label{modifiedpoles}
\end{equation}
with $\sigma = 4$,
and this is the formula used for the numerical results presented here.
Since the spacing still scales with $\sqrt n$, our theoretical
guarantee of convergence still applies.

{\em 2b. Fix\/ $N_2+1 = O(n)$ monomials $1,\dots, (z-z_*)^{N_2}$
and set $N=2 N_1+2 N_2+1$.} For simplicity, we usually take
$N_2=\lceil n/2 \rceil$.  From the point of view of asymptotic
theory one could get away with $N_2 = O(\sqrt n \kern 1pt)$, but
little speed would be gained from such fine tuning and it would
risk delaying the asymptotic convergence in cases dominated by the
smooth part of the problem rather than the corner singularities.
The explanation of the formula for $N$ is that this is the number
of real degrees of freedom in (\ref{form}), since the imaginary
part of the constant term can be taken to be zero.

{\em 2c. Choose $M\approx 3N$ sample points on boundary,  also
clustered near corners.} An essential feature of our algorithm is
that sample points on the boundary are clustered exponentially near
the corners, like the poles $\{z_k\}$.  For a Laplace solution to
accuracy $10^{-8}$, for example, there will be poles at distances
on the order $10^{-8}$ from the corners, and tracking their
effect will require resolution of the boundary on the same scale.
Uniform sampling of the boundary is out of the question, since
it would lead to matrices with on the order of $10^8$ rows.
Indeed, it would not even be a good idea mathematically, since
our algorithm relies on the discrete least-squares norm over the
boundary approximating the continuous supremum norm, and this
property would fail if nearly all of the sample points lay far
from the corner singularities.

Fortunately, experiments show that successful results
do not depend on fine details of the boundary sampling scheme;
and in any case, step (3) of the algorithm gives the chance to
confirm the accuracy of a computed approximation.  The following
is the scheme we have used for our computations.  Each pole $z_k$
is associated with two free parameters of the approximation,
the coefficients multiplying the real and imaginary parts of
$1/(z-z_k)$.  For each $z_k$ we introduce six sample points
on the boundary $\Gamma$.  If $w_j$ is the corner near $z_k$,
define $\delta = |z_k-w_j|$.  We take the six sample points to
be the points at distances $\delta/3$, $2\kern .5pt \delta/3$,
and $\delta$ from $w_j$ along the boundary arcs connecting
$w_j$ with $w_{j-1}$ and $w_{j+1}$.  If any of these points
lie beyond the end of the boundary arc, they are of course not
included.\footnote{This kind of boundary sampling might have
simplified the computations of~\cite{hochman}.}

{\em 2d. Evaluate at sample points to obtain $M\times N$
matrix $A$ and $M$-vector $b$.} Both $A$ and $b$ are real,
with $b$ corresponding to samples of the real function $h(z)$
and the columns of $A$ corresponding to samples of the real and
imaginary parts of the poles and monomials.  In a typical small-
or medium-scale problem, $A$ might have 200--2000 columns and
about three times as many rows.

{\em 2e. Solve the least-squares problem $Ax \approx b$ for the
coefficient vector $x$.} To solve the system, first we rescale the
columns to give them equal 2-norms.  (Without rescaling, columns
$1/(z-z_k)$ corresponding to values $z_k$ very close to a corner
would have some very large entries.)  We then apply standard
methods of numerical linear algebra to solve the least-squares
problem: in MATLAB, a call to the backslash operator, which
applies QR factorization with column pivoting.  The matrices $A$
are often highly ill-conditioned, a phenomenon investigated by
various authors, related to the general notions of over-complete
bases and frames~\cite{bb,hipt,huyb,kit}; see in particular section
3 of~\cite{huyb}.  We make some remarks about the ill-conditioning
of $A$ at the end of Section~\ref{numer} and intend to investigate
this matter further in the future.

{\em 2f. Exit loop if\/ $\|Ax-b\|_\infty < \varepsilon$ or if\/
$N$ is too large or the error is growing.} For problems that are
not too difficult, convergence to the tolerance usually happens
quickly.  If $\varepsilon$ is very small, however, failures become
more likely, and our code includes termination conditions for
such cases.  In particular, our current code becomes unreliable
if $\varepsilon \ll 10^{-10}$, though we hope to improve this
figure with further investigation.

{\em 3. Confirm accuracy by checking the error on
a finer boundary mesh.} One of the appealing features of this
numerical method is that the solution it delivers is exactly a
harmonic function (apart from rounding errors), which implies by
the maximum principle that the error in the interior is bounded
by its maximum on the boundary.  Thus one can get a virtually
guaranteed a posteriori error bound by sampling the boundary error
finely --- say, twice or four times as finely as the least-squares
grid.  Using methods of validated numerics, it should also be
possible to generate a truly guaranteed bound~\cite{tucker},
though we have not pursued this idea.

{\em 4. Construct a function to evaluate $r(z)$ based on the
computed coefficients $x$.}  Once the coefficient vector $x$
for the numerical solution $r(z)$ to the Laplace problem has been
found, many computer languages make it possible to construct an
object~{\tt r} that evaluates the solution at a single point or
a vector or matrix of points with a syntax as simple as \kern
1pt {\tt r(z)}.  In MATLAB this can be done by creating {\tt r}
as an anonymous function.

\section{\label{numer}Numerical examples}
We now present two pages of examples to illustrate the Laplace
solver in its basic mode of operation: a Dirichlet problem on a
polygon with continuous boundary data.  The next section comments on
variants such as Neumann or discontinuous boundary data.
Over the course of the work leading to this article, we have
solved thousands of Laplace problems.

\begin{figure}[ph]
\vskip .15in
\begin{center}\includegraphics[scale=.85]{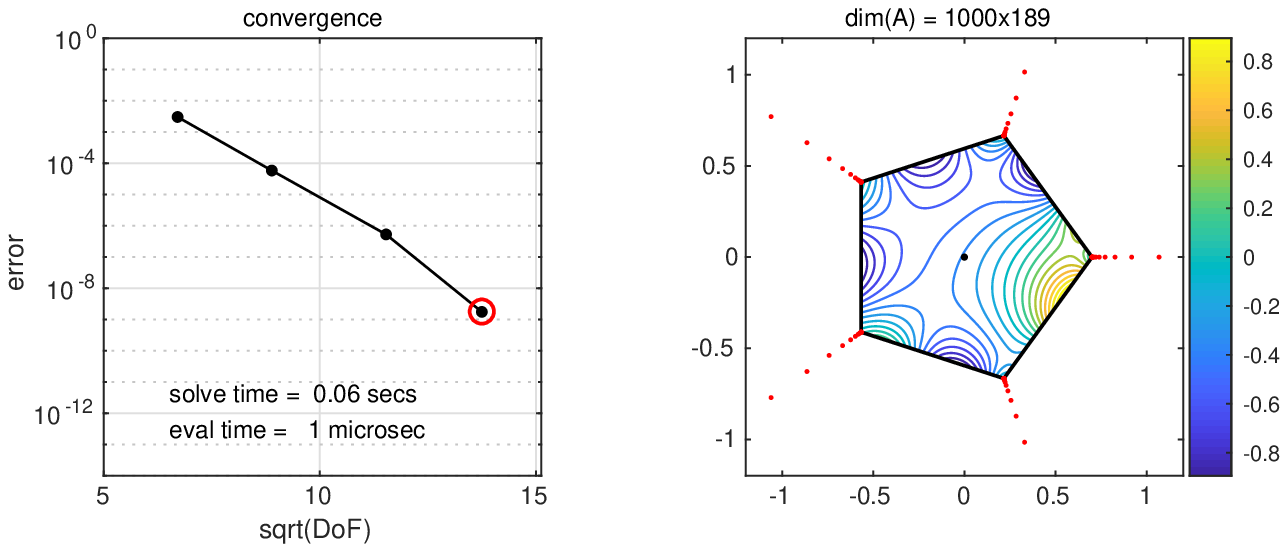}\end{center}
\vskip .15in
\begin{center}\includegraphics[scale=.85]{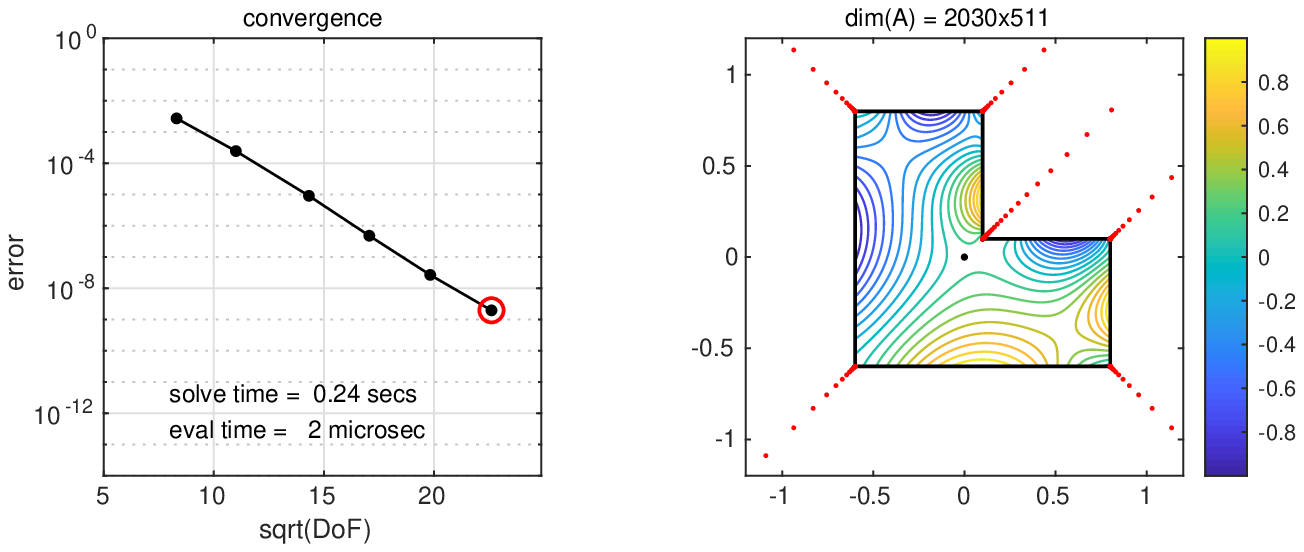}\end{center}
\vskip .15in
\begin{center}\includegraphics[scale=.85]{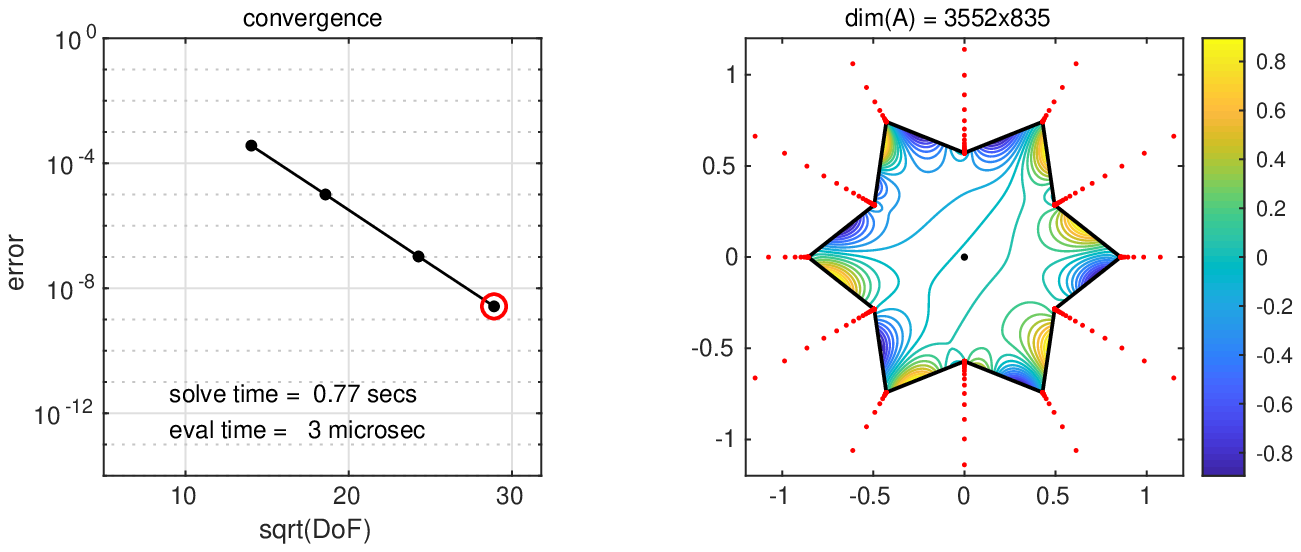}\end{center}
\vskip .15in
\begin{center}\includegraphics[scale=.85]{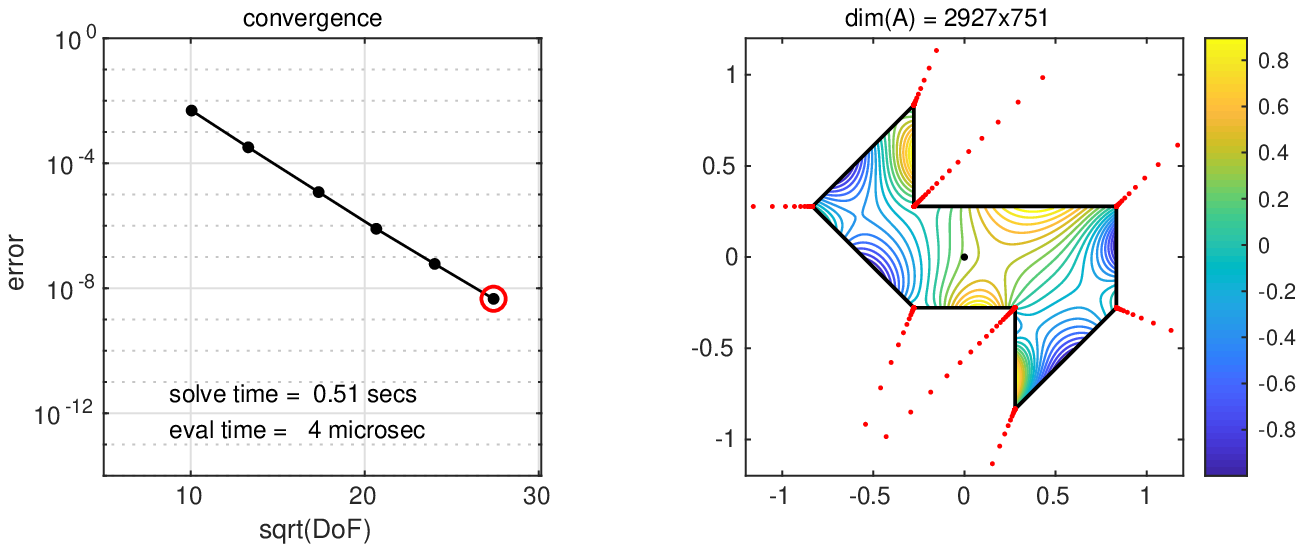}\end{center}
\caption{Four examples of Laplace Dirichlet problems on
polygons: a pentagon, an L-shape, a snowflake, and an isospectral
drum.  The boundary data are analytic functions on each side
with continuity but no further smoothness at the corners.  The
error curves show root-exponential convergence, and the
red circles confirm that the final black dot is truly an upper
bound on the error.}
\label{first4}
\end{figure}

Each example consists of a pair of figures: a convergence curve
on the left and a contour plot of the solution on the right.
Consider the first example of Figure~\ref{first4}, a regular
pentagon.   The convergence curve shows the supremum-norm error
measured over the least-squares grid on the the boundary on a
log scale against $\sqrt N$, where $N$ is the total number of
degrees of freedom (the column dimension of $A$).  The straight
shape of the curve confirms root-exponential convergence down
to the specified tolerance, $\varepsilon = 10^{-8}$.  For this
example the tolerance has been reached with $N=189$ after $0.06$
seconds of desktop time.  The code also does a thousand point
evaluations of the resulting solution $r(z)$ to estimate the time
per point required, which comes out as $1$ microsecond.
Thus one could evaluate the computed solution at ten thousand
points, all with 8-digit accuracy, in a hundredth of a second.

The red circle at the end of the convergence curve corresponds to
a measurement of the boundary error on a finer grid than is
used for the least-squares solve (twice as fine), following
step 3 of the algorithm as described in the last section.  The fact
that it matches the final black dot confirms that the least-squares
grid has been fine enough to serve effectively as a continuum.

The image on the right shows the computed solution.  In this set
of examples, the boundary function $h$ is taken to be a different
analytic function on each side taking the value zero at the corners
to ensure continuity.  (The boundary functions are smooth random
functions of the kind produced
by the Chebfun {\tt randnfun} command~\cite{smooth}, but this is just for
convenience, not important to our algorithm.)  This figure also
shows the poles exponentially clustered outside each corner as
red dots.  The distances
go down to the same order of magnitude as the accuracy, $10^{-8}$, so most
of the poles are indistinguishable in the figure.  There is also
a black dot at the origin indicating the expansion point $z_*$ of the
polynomial part of (\ref{form}).  The heading reports the
dimension of the final matrix used, whose least-squares solution
dominates the computing cost.

The remaining examples of Figure~\ref{first4} show an
L-shape, a snowflake, and an isospectral drum.
Root-exponential convergence is evident in each case.
Each problem is solved in less than a second, and the
solutions obtained can be evaluated in
less than 5 microseconds per point.

\begin{figure}[ph]
\vskip .15in
\begin{center}\includegraphics[scale=.85]{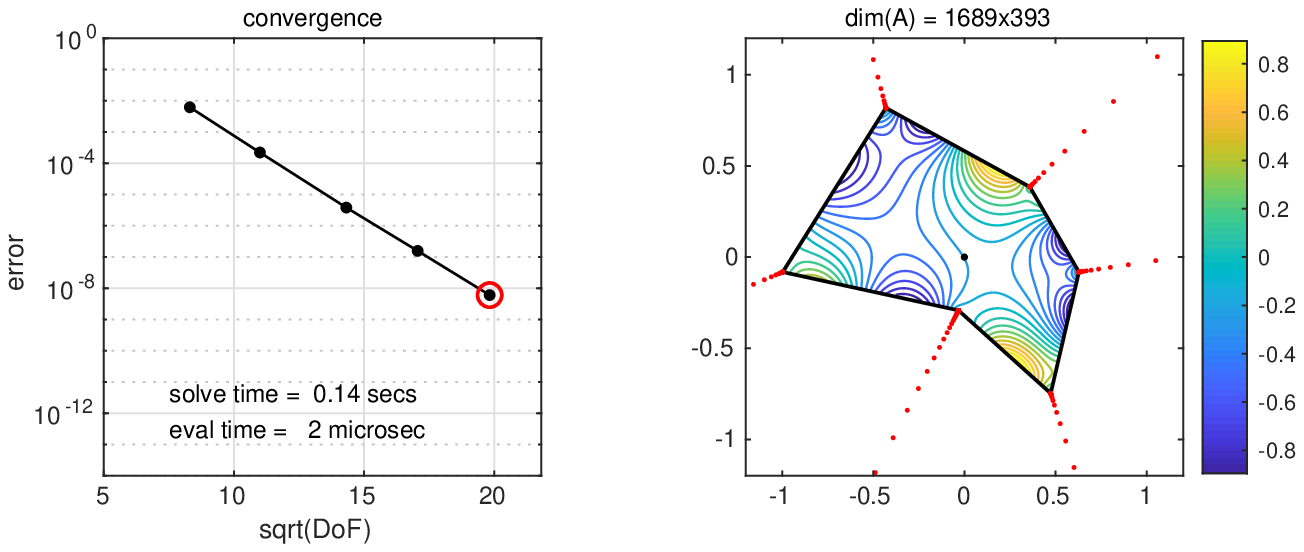}\end{center}
\vskip .15in
\begin{center}\includegraphics[scale=.85]{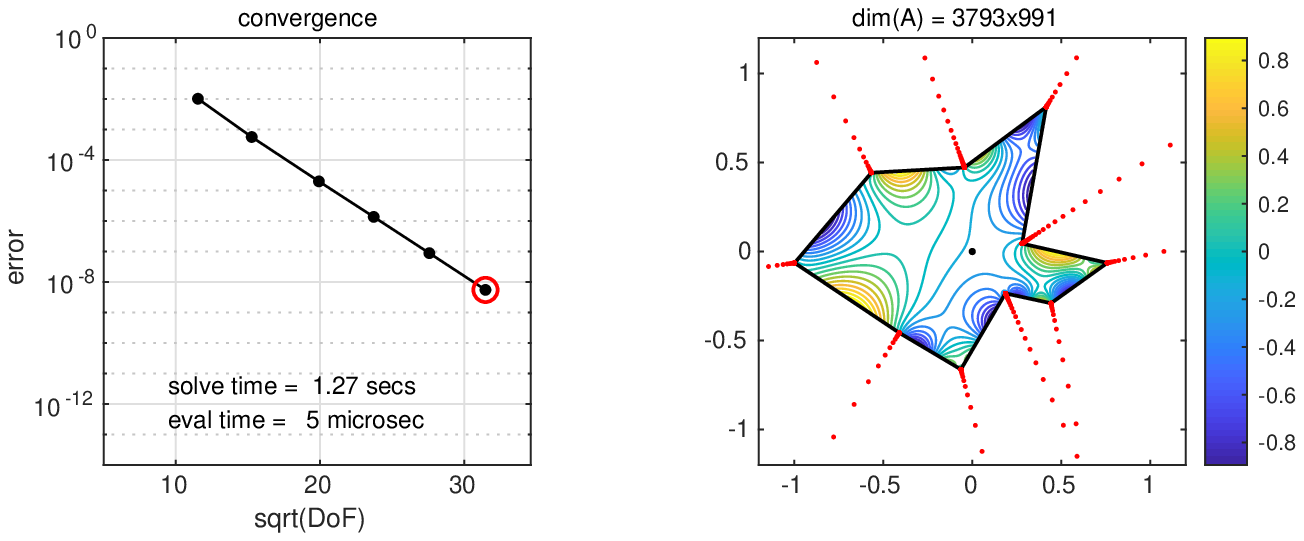}\end{center}
\vskip .15in
\begin{center}\includegraphics[scale=.85]{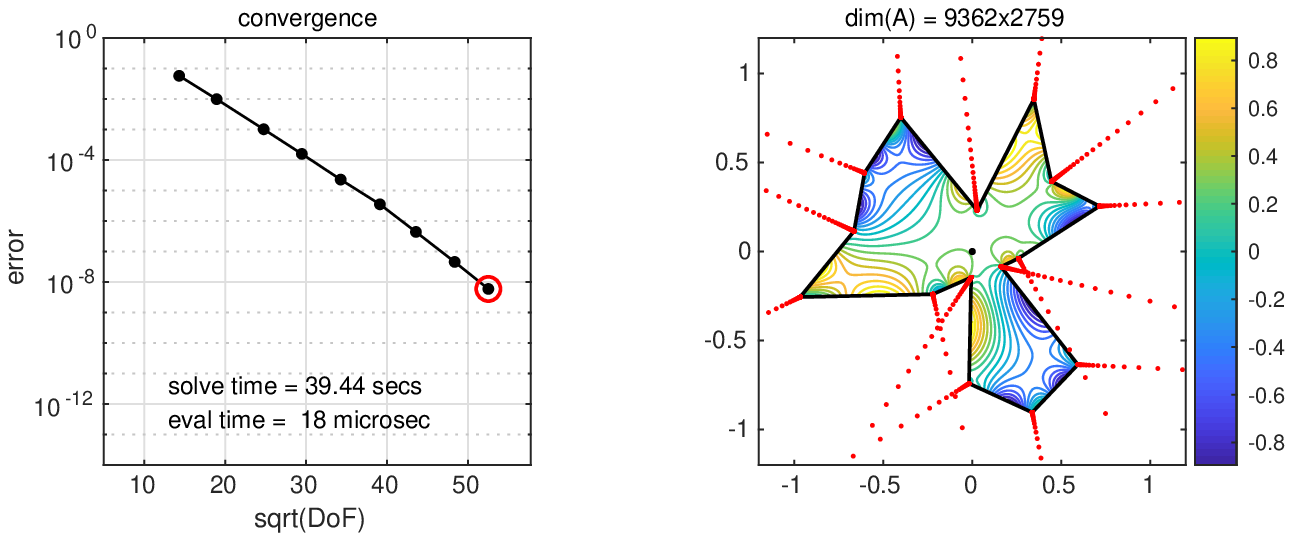}\end{center}
\vskip .15in
\begin{center}\includegraphics[scale=.85]{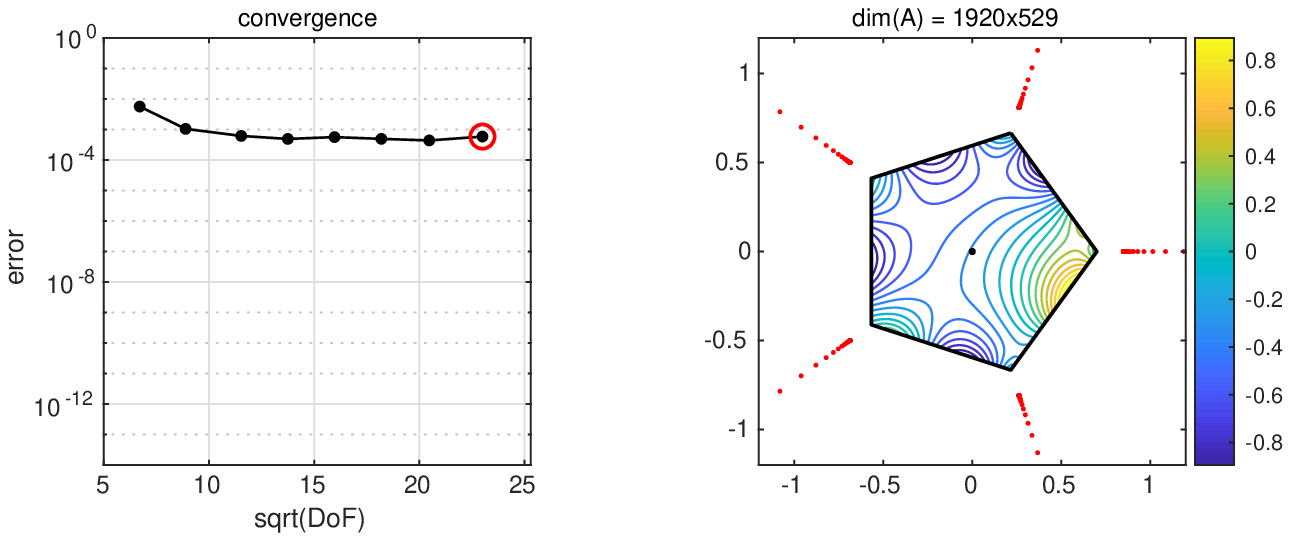}\end{center}
\caption{Four more examples.  The first three involve random polygons
with $6$, $10$, and $15$ sides.
With $\varepsilon=10^{-4}$ instead of $10^{-8}$, the
solve times of the latter two improve to
$0.18$ and $2.3$ seconds.  The final image shows that
root-exponential convergence fails if the poles are
displaced a small distance away from the corners.}
\label{second4}
\end{figure}

The second page of examples, Figure~\ref{second4}, first shows random
polygons with 6, 10, and 15 vertices.  The vertices of each
polygon lie at uniformly spaced angles about a central point, with the
radii taken as independent samples from a uniform distribution.
The figures show successful
computations, but the times slow down cubically when there
are many vertices.  For the 15-gon, one would do
better to use the advanced integral equations methods
of Serkh and Rokhlin~\cite{ser} or Helsing and
Ojala~\cite{helsing,helsingo,ojala}, who solve problems with thousand
of corners with their method of recursive compressed inverse
preconditioning applied on graded meshes.  On the other hand,
loosening the tolerance in Figure~\ref{second4} 
speeds up the computations greatly.  With $\varepsilon
= 10^{-4}$ instead of $10^{-8}$, the solve times for the last 
two examples improve from $1.3$ and $39$ seconds to
$0.18$ and $2.3$ seconds. 

\begin{figure}
\begin{center}
\includegraphics[scale=.85]{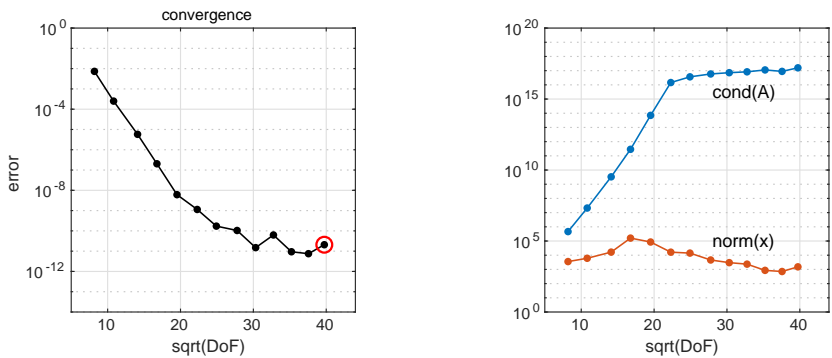}
\end{center}
\caption{\label{illcond} The second example of
Figure\/~$\ref{second4}$ carried to $1579$ degrees
of freedom.  The accuracy
stagnates at around $10^{-11}$, but the curves
on the right show that the matrix $A$ (after column
scaling) reaches a condition
number of on the order of $10^{16}$ long
before this point.}
\end{figure}

The last image of Figure~\ref{second4} attempts a computation with
all the poles shifted a distance $0.15$ away from their corners.
The fast convergence is now lost completely, confirming the
decisive importance of Newman's phenomenon and
of our analysis in Section~2.

These figures show the speed and accuracy of the rational
approximation algorithm.  One could look deeper, however, to assess
its numerical properties and limitations.  In particular, it was
mentioned in the last section that our current implementation
cannot reliably compute solutions with accuracy better than
$10^{-10}$.  It is natural to ask, where does this limitation
come from, and could the situation be improved?  At present, our
understanding is not sufficient to attempt a proper discussion of
this question.  Certainly a relevant fact is that our least-squares
problems lie in the regime where the matrices are hugely
ill-conditioned and yet reasonable fits are obtained anyway.
To illustrate this, Figure~\ref{illcond} shows the second example
of Figure~\ref{second4}, but carried to higher accuracy.  Instead
of six data points, there are now 13, and one sees stagnation
around accuracy $10^{-11}$.  The plot on the right shows that long
before this point, the matrix $A$ has reached a condition number
on the scale of the reciprocal of machine precision, though the
norms of the coefficient vectors remain under control.  This is a
familiar effect with overcomplete bases, investigated by various
authors over the years~\cite{bb,hipt,huyb,kit}.  We believe that further
analysis will provide an understanding of how these effects play
out in the present context and suggest modifications to bring
our results closer to machine precision.  At the same time,
it is worth emphasizing that without such analysis,
we still get fast solutions to excellent accuracy.

\section{Variants}

We now mention a number of variations on the
basic problem considered so far.

\begin{figure}
\begin{center}
\begin{center}\includegraphics[scale=.85]{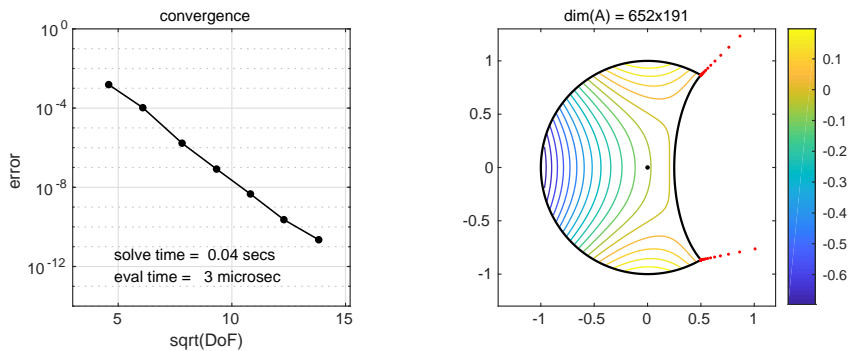}\end{center}
\end{center}
\caption{\label{variants}Illustration of a domain with
curved boundaries.  The strings of poles have been placed askew
to illustrate that precise positioning is unimportant, so
long as the clustering is exponential.}
\end{figure}

{\em Curved boundaries.}  There is nothing new to be done here.
The method we have described works with curved boundaries, the only
difficulty being that the programming is less straightforward
since the definition of the geometry involves more than just
line segments.  
Fig.~\ref{variants} shows the solution to ten-digit accuracy of a
Dirichlet problem of this kind on a domain bounded a circular
arc and an elliptical arc.
The strings of poles have been placed along lines oriented
at arbitrary angles to illustrate that using precisely the external
bisector is not important.

{\em Neumann boundary conditions.}  Little changes here; the
boundary matching now involves derivatives as well as function
values.  For a rigorous a posteriori error bound, one will
have to make sure to use a formulation of the maximum principle that
applies with the given boundary conditions.

{\em Transmission problems.}  Some problems involve transmission
of a signal from one domain to another with a matching condition
at the interface, often involving derivatives.  In this case we
would use distinct rational representations in the different regions.

{\em Discontinuous boundary conditions.}
In applications, Laplace problems very often have
discontinuous boundary conditions; a prototypical example would be
a square with boundary values $0$ on three
sides and $1$ on the fourth.  Mathematically, this is
no problem at all, and the same applies to our numerical algorithm.
The change that must be introduced is that convergence of
the error to zero in the supremum norm is no longer possible.
Instead, we find it convenient to guide the computation
by a supremum norm weighted by distance to the nearest corner.
Again one must be careful in applying the maximum
principle.

{\em Long domains.}  The ``Runge part'' of
(\ref{form}) consists of a monomial expansion about a
point $z_*$.  If $N_2$ is large and $\Omega$ is highly elongated,
one is likely to run into trouble because this is an
exponentially ill-conditioned basis.  Approximation theory
suggests turning to the basis of the
Faber polynomials for $\Omega$~\cite{gaier},
but this would introduce an auxiliary computation as difficult
as the original Laplace problem.  We believe that simpler methods
would be effective in many cases, such as a switch to a suitably
scaled and shifted Chebyshev basis if $\Omega$ is elongated mainly
in one direction.
Impressive solutions for such problems based on integral
equations can be found in~\cite{ojala}.

{\em Poisson equation.} For a problem $\Delta u = f$ with boundary
conditions $u(z) = h(z)$ for $z\in \Gamma$, as is very well
known, one can first find any function
$v$ such that $\Delta v = f$ in $\Omega$, with arbitrary boundary values.
The correction $w = u-v$ needed to find $u$ is now the solution of
the Laplace problem $\Delta w = 0$ with boundary conditions $w(z)
= h(z) - v(z)$.  This Laplace problem will have
corner singularities as usual.

{\em Helmholtz equation.}  Here, it appears that rational
functions can be generalized
to series involving Hankel functions, though
the mathematics brings significant challenges.  A preliminary announcement
appeared in~\cite{pnas}, and development is underway.

{\em Multiply connected domains.} In a multiply connected domain,
poles will need to be placed within the holes as well as in the
exterior (and also logarithmic terms associated with
the Runge part of~(\ref{form})~\cite{axler,series}).
In the context of rational functions,
as observed by Runge himself~\cite{gaier,runge}, this is not a very
significant change.  Of course, as always, one will have to
make sure that no poles fall within~$\Omega$.
Illustrations of series methods for Laplace
problems in multiply connected smooth regions, without exponential
clustering of poles, can be found in~\cite{series}.

{\em Slits.}  A domain with a slit, that is, a corner with
interior half-angle $\pi$, appears problematic since one
can hardly place poles on
the slit.  However, local transformations can get around such
a problem.  If $\Omega$ has a boundary slit coming in to
$z=0$ along the negative real axis, for example, then the
local change of variables $\zeta = z^{1/2}$ transforms the
slit to a segment of the imaginary axis.  One can now introduce
poles $\zeta_k$ in the $\zeta$-plane to get terms for a Newman-style
expansion, with $1/(\zeta-\zeta_k)$ corresponding to
$1/(z^{1/2} -\zeta_k)$.
Such transformations are just local, changing only the 
nature of some of the terms appearing in (\ref{form}); the
least-squares problem is still solved on
the boundary $\Gamma$ in the $z$-plane.

{\em Faster than root-exponential convergence.}  
As mentioned at the beginning of Section~\ref{thms}, faster than
root-exponential convergence is not possible for Laplace problems
in regions with corners, if accuracy is measured
by the maximum norm over $\Omega$.
However, as pointed out in Myth~3
of~\cite{myths} and Chapter~16 of~\cite{atap},
best approximations may not be optimal!  The point
here is that the insistence on global supremum norm error bounds may
have terribly adverse effects on an approximation throughout almost
all of a domain.    Specifically, the discussion of discontinuous
boundary conditions above suggests a possibility that might be very
appealing to users.  Even for a problem with continuous boundary
data, one could measure convergence by a supremum norm weighted
by distance to the nearest corner.  This could make it possible to
accelerate convergence of our algorithm
from root-exponential to nearly exponential,
yet still guarantee errors bounded by a tolerance divided by the
distance to the nearest corner.  We plan to include an option
along these lines in our software.

\section{Discussion}
The numerical solution of PDE\kern .5pt s has been brought to a
very highly developed state by generations of mathematicians and
engineers since the 1950\kern .5pt s.  The Laplace equation in
2\kern .5pt D is as basic a problem in this area as one could ask
for, and many methods can handle it.  Getting high accuracy in the
presence of corner singularities, however, still requires care,
as we found from the responses to a challenge problem involving
the L-shaped domain that we posed to the NA Digest email discussion
group in November, 2018~\cite{pnas,NA}.  The most common approach
to such a problem would be to use the finite element method
(FEM), for which there is widely-distributed freely-available
software such as deal.II, FEniCS, Firedrake, IFISS, PLTMG, and
XLiFE++~\cite{dealII,pltmg,ifiss,fenics,firedrake,xlife}.  However,
it is not straightforward to get, say, 8 digits of accuracy by such
methods.\footnote{Just one respondent to our posting solved the
problem to more than 8 digits of accuracy by FEM, Zo\"is Moitier
of the University of Rennes~1, who used XLiFE++~\cite{xlife} with
finite elements of order 15.} For this, it may be preferable to
use the more specialized tools of $hp$-adaptive FEM~\cite{schwab}.
Such methods are powerful, but we believe they cannot compete
for speed and simplicity for these simple planar problems, and
in particular, we are unaware of any FEM software that can match
the performance shown in Figures~\ref{first4} and~\ref{second4}.

The other highly developed approach for Laplace problems is
boundary integral equations (BIE), more restricted than FEM
but very powerful when applicable~\cite{rokhlin}.  Whereas FEM
constructs a two-dimensional representation of the solution of
a PDE, BIE represents it via an integral over a one-dimensional
density function on the boundary.  This speeds up the computations
greatly, and in addition, one can often avoid the ill-conditioned
matrices associated with FEM, especially through
formulations involving Fredholm integral equations of the second kind.
Evaluations of a BIE solution
require numerical quadrature, which poses challenges near the
corners and also near the boundary away from the corners, but
experts have developed quadrature methods that are fast and
accurate~\cite{brs,qbx,ser}.  We have already mentioned the method
of Helsing and Ojala with its ability to solve problems with
thousands of corners~\cite{helsing,helsingo,ojala}.
There is no doubt that integral equations constitute the most
powerful tool currently available for 2D Laplace problems.

Our method, which we think of as constructing a zero-dimensional
representation of the solution, belongs to the category
of techniques known as the Method of Fundamental Solutions
(MFS)~\cite{bb,bog,fair,kat,kit,liu}.  This refers to any
method in which a function is found satisfying given boundary
conditions by expansion in a series of free space solutions with
point singularities.  Sometimes the points are fixed in advance,
making the calculation linear as in our method, and in other
cases they are selected adaptively.  Such ideas go back at least
to Kupradze in the 1960\kern .3pt s~\cite{kup,kupbook} and are
a special case of so-called Trefftz methods~\cite{hipt,hipt4}.
One difference of what we are proposing here from the usual is
that MFS expansions usually involve monopoles (= logarithmic point
charges) rather than dipoles (= simple poles of a meromorphic
function).  A second difference is that the MFS literature does
not make the connection with rational approximation theory, and in
particular, does not normally use exponentially clustered poles
to achieve root-exponential convergence near singularities,
though steps in this direction can be found in~\cite{liu}.
Another related method is that of Hochman, et al.~\cite{hochman},
involving rational functions with a nonlinear iteration.

It is interesting to consider the conceptual links between
boundary integrals and MFS methods such as our rational functions.
MFS methods seek to expand the solution in a finite series of
point singularities lying outside the problem domain.  Boundary
integral equations, by contrast, seek a continuum distribution of
charge singularities (single-layer for monopoles, double-layer
for dipoles) precisely along the boundary.  This formulation on
the boundary has the advantage that it is connected with highly
developed and powerful mathematical theories of integral operators,
while point charge approaches are less fully understood.  On the
other hand the BIE approach has the flavor of interpolation rather
than least-squares, requiring exact determination of a precisely
determined density function, quite a contrast to the great
simplicity and convenience of least-squares problems of the point
charge methods, which take advantage of ``overcomplete bases.''
For example, Walsh and his student Curtiss\footnote{Donald Newman
was also a student of Walsh's, receiving his
PhD from Harvard at age 22.} were interested
in solving Laplace problems by interpolatory boundary matching,
and ended up giving a good deal of their effort to the interpolation
rather than the approximation theory~\cite{curtiss,walsh29}.
The much simpler idea of using least-squares boundary matching
goes back to Moler and undoubtedly others in the 1960\kern .3pt
s if not before~\cite{moler,series}, but has been surprisingly slow to find a
place in mainstream numerical PDE\kern .5pt s.

As commented in~\cite{pnas}, we see a historical dimension
in the relative lack of attention given to global methods of
a spectral flavor for solving PDE\kern .5pt s, such as we have
proposed here, that take advantage of the analyticity of problems
arising in applications.  In the 19th century, mathematicians'
focus was mainly on analytic functions, possibly with branch
points.  In the 20th century, however, interest in that kind of
function theory largely gave way to interest in real analysis,
with great attention being given to detailed study of questions of
regularity, that is, of precise degrees of smoothness of functions
and their consequences.  This trend is strikingly on display in
the standard work on elliptic PDE\kern .5pt s in domains with
corners, by Grisvard~\cite{grisvard}, which begins with 80 pages
on Sobolev spaces and proceeds to develop an extensive theory in
that setting.  The numerical PDE community followed this trend,
with the prevailing discourse again becoming analysis in terms
of Sobolev spaces~\cite{schwab}.  Although in principle one can
analyze even analytic problems with these tools, in practice they
bring a bias toward discretizations tuned to problems of limited
smoothness.\footnote{For example, in the FEM literature as in
Grisvard's book one encounters the view that a PDE problem with
a salient corner is nonsingular whereas the same problem with
a reentrant corner is singular.  Such a distinction would not
make sense from the point of view of classical function theory.}
Yet many problems of interest in applications are analytic or
piecewise analytic.

Turning to future prospects, we note that the 2\kern .3pt D
Laplace equation is very special.  We have begun exploration of
generalizations to the Helmholtz equation~\cite{pnas}, and other
problems such as the biharmonic equation can also be investigated.
The extension to 3D is a bigger issue to be considered.  Here,
as in 2\kern .3 pt D, there are MFS papers going back many
years.  The question is whether a new more focussed connection
with approximation theory near singularities can bring a new
efficiency, as we have shown here for 2\kern .3pt D problems,
and we are open-minded on this question.

In this paper we have introduced a new approach to numerical
solution of certain PDE\kern .5pt s, globally accurate and
extremely simple, requiring no preliminary analysis of corner
singularities.  The method is very young, and there are innumerable
questions to be answered and details to be improved in further
research.  Since the method exploits the same mathematics that
makes lightning strike trees and buildings at sharp edges, we
like to think of it as a ``lightning Laplace solver.''

\section*{Acknowledgments}
We are grateful for advice to Alex Barnett, Timo Betcke, Leslie
Greengard, Dave Hewett, Daan Huybrechs,
Yuji Nakatsukasa, Vladimir Rokhlin,
Kirill Serkh, Andr\'e Weideman, and an anonymous referee.
We have also benefited greatly from Bernhard Beckermann's
extensive knowledge of the rational approximation literature.

\end{document}